\documentclass[11pt,halfparskip]{scrartcl}

\usepackage{amsfonts}
\usepackage{amsmath, amsthm}
\usepackage{nicefrac}

\usepackage{graphicx}
\usepackage{subfigure}

\usepackage{algpseudocode}
\usepackage{algorithm}

\usepackage{hyperref}
\usepackage[T1]{fontenc}
\usepackage{textcomp}

\newtheorem{thm}{Theorem}

\newtheorem{cor}{Corollary}

\newtheorem{lem}{Lemma}

\newtheorem{defi}{Definition}                 

\newtheorem{exa}{Example}

\newcommand{\argmax}{{\mathrm{argmax}}\,}

\newcommand{\inte}{{\mathrm{int}}\,}
\newcommand{\relint}{{\mathrm{relint}}\,}

\newcommand{\conv}{{\mathrm{conv}}\,}

\newcommand{\prox}{{\mathrm{proj}_X}\,}

\newcommand{\R}{{\mathbb R}}

\newcommand{\N}{{\mathbb N}}

\newcommand{\Q}{{\mathbb Q}}
\newcommand{\Z}{{\mathbb Z}}

\mathchardef\ordinarycolon\mathcode`\:
\mathcode`\:=\string"8000
\begingroup \catcode`\:=\active
  \gdef:{\mathrel{\mathop\ordinarycolon}}
\endgroup

\title{$k$-disjunctive cuts and a finite cutting plane algorithm for
  general mixed integer linear programs}

\author{Markus J\"org \\
  Technische Universit\"at M\"unchen, Zentrum Mathematik \\
  Boltzmannstra{\ss}e 3, 85747 Garching bei M\"unchen, Germany \\
  joerg@ma.tum.de }

\date{July 26, 2007}

\begin{document}
\maketitle
\begin{abstract}
  \textbf{Abstract:} In this paper we give a generalization of the
  well known split cuts of Cook, Kannan and Schrijver \cite{CKS90} to
  cuts which are based on multi-term disjunctions. They will be called
  $k$-disjunctive cuts. The starting point is the question what kind
  of cuts is needed for a finite cutting plane algorithm for general
  mixed integer programs. We will deal with this question in detail
  and derive cutting planes based on $k$-disjunctions related to a
  given cut vector. Finally we will show how a finite cutting plane
  algorithm can be established using these cuts in combination with
  Gomory mixed integer cuts.
\end{abstract}

\section{Introduction}\label{sec:intro}
In this paper we will deal with cutting planes and related algorithms
for general mixed integer linear programs (MILP). As most of the
results will be derived by geometric arguments we focus on programs
that are given by inequality constraints, i.e.
\begin{equation}\label{eq:MILP}
\begin{array}{ll}
  \max & cx + hy \\
       & Ax + Gy \leq b \\
       & x \in \Z^p
\end{array}
\end{equation}
where the input data are the matrices $ A \in \Q^{m \times p}, G \in
\Q^{m \times q}$, the column vector $b \in \Q^m$ and the row vectors $
c \in \Q^p, h \in \Q^q$.  Moreover we denote by the polyhedra $P = \{
(x,y) : Ax + Gy \leq b\} \subset \R^{p+q}$ and $ P_I = \conv( \{(x,y)
\in P: x \in \Z^p\}) \subset \R^{p+q} $ the feasible domains of the LP
relaxation and the (mixed) integer hull of a given MILP, respectively.
We call a MILP bounded, if the polyhedron $P$ is bounded. We will also
need the projection $\prox(P) := \{ x \in \R^p : \exists y \in \R^q :
(x,y) \in P\}$ of the polyhedron $P$ on the space of the integer
variables.

By a cutting plane for $P$ we understand an inequality $ \alpha x +
\beta y \leq \gamma$ with row vectors $ \alpha \in \Q^p, \beta \in
\Q^q$ which is valid for $P_I$ but not for $P$. Using cutting planes
gives a simple idea of how to solve a general MILP: Solve the LP
relaxation of the MILP. If the optimal solution is feasible, i.e.
satisfies the integrality constraint, an optimal solution is found.
Otherwise find a valid cutting plane that cuts off the current
solution and repeat. But unlike the pure integer case no finite exact
cutting plane algorithm is known for general MILP. Therefor we remark
that most cutting planes for general MILP such as e.g.  Gomory mixed
integer cuts \cite{Gom63} or mixed integer rounding cuts \cite{NW90}
are special cases of or equivalent to split cuts \cite{CKS90}. This
fact and more detailed relations between these and other cuts are
stated in \cite{CL01}.  Here a split cut is defined as a cutting plane
$ \alpha x + \beta y \leq \gamma$ for $P$ with the additional property
that there exists $ d \in \Z^p, \delta \in \Z$ such that $ \alpha x +
\beta y \leq \gamma$ is valid for all $ (x,y) \in P$ which satisfy the
split disjunction $ dx \leq \delta $ or $ dx \geq \delta + 1$. So
split cuts are defined not constructively but by a property, only. Now
one can see in the following 'classical' example of Cook, Kannan and
Schrijver \cite{CKS90} that split cuts are not sufficient for solving
a general MILP in finite time.
\begin{exa}\label{exa:cks}
  The MILP 
  \begin{eqnarray*}
    && \max y \\
    && -x_1 + y \leq 0 \\
    && -x_2 + y \leq 0 \\
    &&  x_1 + x_2 + y \leq 2 \\
    && x_1, x_2 \in \Z
  \end{eqnarray*}
  has the optimal objective function value 0 but the problem cannot be
  solved by any algorithm that uses split cuts, only. A proof of this
  statement in a more general context is given in
  \autoref{lem:expneed}.
\end{exa}
On the other hand, as positive results in the context of cutting plane
algorithms for MILP we can only give the following two special cases:
For mixed 0-1 programs split cuts are sufficient for generating the
integer hull $P_I$ of a given polyhedron $P$. See e.g. \cite{NW90} in
the context of mixed integer rounding cuts or \cite{BCC93} in the more
recent representation of lift-and-project cuts.  For general MILP,
there only exists a finite approximation algorithm of Owen and
Mehrotra \cite{SM01} which finds a feasible $\epsilon$-optimal
solution and uses simple split cuts, that means split cuts to
disjunctions $x_i \leq \delta \vee x_i \geq \delta + 1$.

So as split cuts fail in the design of a finite cutting plane
algorithm for general MILP we want to generalize this approach to cuts
that are based on multi-term disjunctions. Therefor we start in
\autoref{sec:kdis} with the introduction of $k$-disjunctive cuts and
some of its basic properties. Afterwards we look at the approximation
properties of the $k$-disjunctive closures and deal with the question
what kind of cuts is needed for an exact finite cutting plane
algorithm both in general and in special cases. Finally we derive a
$k$-disjunctive cut according to a given cut vector. In
\autoref{sec:algo} we turn to algorithmic aspects and give a way of
how a finite cutting plane algorithm for general MILP can be designed
using $k$-disjunctive cuts in connection with the well known mixed
integer Gomory cuts. Finally we will discuss the algorithm and give
some interpretations.

\subsection{Preliminaries}\label{ssec:prel}
Here we repeat two basic results that we will need during this paper.
The first one deals with the computation of the projection $\prox(P)$,
the second one with the convergence of the mixed integer Gomory
algorithm in a special case.
\begin{lem}\label{lem:projection}
  Let a polyhedron $P = \{ (x,y) : Ax + Gy \leq b\}$ be given. Then
  \begin{equation*}
     \prox(P) = \{x \in \R^p:\ v^rAx \leq v^rb, \ \forall r \in R \},
  \end{equation*}
  where $\{v^r\}_{r \in R} $ is the set of extreme rays of the cone $
  Q := \{v \in \R^m: \ G^Tv = 0, v \geq 0\}.$
 \end{lem}
 \begin{proof}
   The statement follows by applying the Farkas Lemma, see e.g.
   \cite{NW88}, I.4.4.
 \end{proof}
 Next we look at the usual mixed integer Gomory algorithm
 \cite{Gom63}. Although the algorithm does in general not even
 converge to the optimum, the special case in which the optimal
 objective function value can be assumed to be integer, e.g. the case
 of $h = 0$, can be solved finitely using the algorithm. In detail we
 have the following
\begin{thm}\label{thm:gomory}
  Let a bounded MILP \eqref{eq:MILP} be given. Then the mixed integer
  Gomory algorithm terminates finitely with an optimal solution or
  detects infeasibility under the following conditions:
  \begin{enumerate}
  \item One uses the lexicographic version of the simplex algorithm
    for solving the LP relaxation.
  \item The optimal objective function value is integral.
  \item A least index rule is used for cut generation, i.e. the
    mixed integer Gomory cut according to the first variable $x_j$,
    that is fractional in the current LP solution, is added to the
    program. Here $x_0$ corresponds to the objective function value.
  \end{enumerate}
\end{thm}
Using the last theorem, it is obvious that we can check in finite time
if there is a feasible point in a polytope with a given (rational)
objective function value, as by scaling it can be always assumed that
the optimal objective function value is integral. This is expressed in
the following
\begin{cor}\label{cor:gomory}
  Let a bounded MILP \eqref{eq:MILP} with the additional constraint $
  cx + hy = \gamma$ be given. Then the mixed integer Gomory algorithm
  terminates finitely with a feasible solution or detects
  infeasibility.
\end{cor}

\section{$k$-disjunctive cuts}\label{sec:kdis}
\subsection{Basic definitions and properties}\label{ssec:basic}
In analogy to the definition of a split cut based on a split
disjunction we now define a $k$-disjunctive cut that is based on a
$k$-disjunction that contains every integral vector.
\begin{defi}\label{defi:kdis}
  Let $k \geq 2$ be a natural number, $d^1,\ldots,d^k \in \Z^p $
  integral vectors and $\delta^1,\ldots,\delta^k \in \Z$. Then we call
  the inequalities $d^1x \leq \delta^1,\ldots,d^kx \leq \delta^k$ a
  $k$-disjunction, if for all $x \in \Z^p$ there is an $i \in
  \{1,\ldots,k\}$ with $ d^ix \leq \delta^i$. In this case we write
  $D(k,d,\delta)$ with $ d = (d^1,\ldots,d^k), \delta =
  (\delta^1,\ldots,\delta^k)$ for the $k$-disjunction.
\end{defi} 
\begin{figure}
  \centering 
  \subfigure[a 3-disjunction \label{sfig:3dis}]
  {\includegraphics[width=45mm]{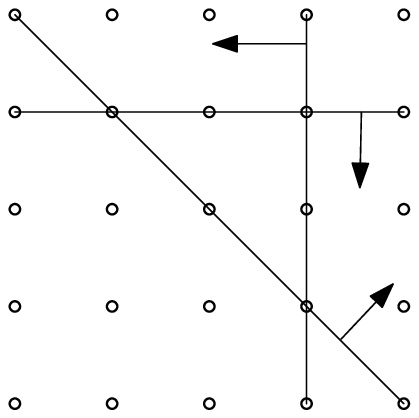}} \hspace*{1.3cm}
  \subfigure[a 4-disjunction \label{sfig:4dis}]
  {\includegraphics[width=45mm]{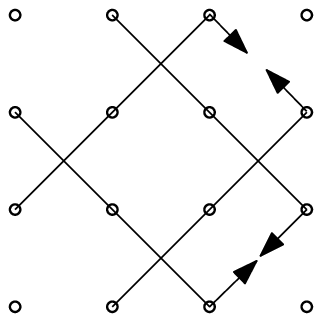}}
  \caption{Examples for $k$-disjunctions in $\R^2$}
  \label{fig:kdis}
\end{figure}
We note that we do not require the vectors $d^i, \delta^i$ to be
different. So every $l$-disjunction is also a $k$-disjunction for $l <
k$. Especially every split disjunction is also a $k$-disjunction.
Moreover every $k$-disjunction is a cover of $\Z^p$ by definition.
\begin{defi}\label{defi:kdiscut}
  Let $P \subset \R^{p+q }$ be a polyhedron and $\alpha x + \beta
  y\leq \gamma$ be a cutting plane. Then $\alpha x + \beta y\leq
  \gamma$ is called a $k$-disjunctive cut for $P$, if there exists a
  $k$-disjunction $D(k,d,\delta)$ with
  \begin{equation*}
    (x,y) \in P: \ \alpha x + \beta y >  \gamma
    \Longrightarrow d^i x > \delta^i, \ \forall i \in \{1,\ldots,k\} .
  \end{equation*}
\end{defi}
Of course every $k$-disjunctive cut for $P$ is valid for $P_I$ by
definition. According to the remark after \autoref{defi:kdis} every
$l$-disjunctive cut is also a $k$-disjunctive cut for $l < k$. So
every split cut is a $k$-disjunctive cut.
\begin{defi}\label{defi:kdisclos}
  Let $P \subset \R^{p+q }$ be a polyhedron. Then the intersection of
  all $k$-disjunctive inequalities is called the $k$-disjunctive
  closure of $P$ and denoted by $P_k^{(1)}$. Analog the $i$-th
  $k$-disjunctive closure $P_k^{(i)}$ of $P$ is defined as the
  $k$-disjunctive closure of $P_k^{(i-1)}$. In the special case of $k
  = 2$ we will also write $P^{(i)}$ instead of  $P_2^{(i)}$.
\end{defi}
We want to remark that it is not evident if the $k$-disjunctive
closure $P_k^{(1)}$ for a given polyhedron $P$ is again a polyhedron
in the case of $k \geq 3$. The both proofs of this property for the
split closure \cite{ACL05} and \cite{CKS90} cannot be applied to the
more general case. However, we will not further deal with this
question, as our results in the following are independent of this
property. We further remark that \autoref{defi:kdiscut} also applies
in a natural way to closed convex sets $P$. Therewith it is guaranteed
that the definition of the $i$-th $k$-disjunctive closure $P_k^{(i)}$
of a polyhedron $P$ is well defined.

A valid cut to a given $k$-disjunction can be computed as intersection
cut to any basis solution of the LP-relaxation that is not contained
in the disjunction according to \cite{Bal71}. In the case of $k = 2$,
Andersen, Cornu\'ejols and Li have shown \cite{ACL05} that
intersection cuts are sufficient to describe all cuts to a given split
disjunction. This result is not true for general $k$-disjunctions.
Here not every valid $k$-disjunctive cut to a given disjunction is
equal to or dominated by an intersection cut. This can be seen in the
following
\begin{exa}
  We look at the polyhedral cone $C \in \R^{2+1}$ with apex
  $(\frac{1}{2},\frac{1}{2},\frac{1}{2})$ that is defined by
  \begin{eqnarray*}
    -x_1 + y  &\leq& 0 \\
    -x_2 + y  &\leq& 0 \\
     x_1 + y  &\leq&  1 \\
     x_2 + y  &\leq& 1.
  \end{eqnarray*}
  Then $ y \leq 0$ is a $4$-disjunctive cut for $C$ to the
  $4$-disjunction $D:= \{x_1 + x_2 \geq 2, x_1 - x_2 \geq 1, -x_1 +
  x_2 \geq 1, -x_1 - x_2 \geq 0 \}$ . The set of all bases is given by
  any three of the above constraints. Computing the four relating
  intersection cuts to the $4$-disjunction $D$ we get that the point
  $(\frac{1}{2},\frac{1}{2},\frac{1}{6})$ is valid for the four cuts
  and so $ y \geq \frac{1}{6}$ has to be satisfied.
\end{exa}
Although the properties of general $k$-disjunctive cuts are more
involved than in the case of split cuts, an investigation of these
cuts is useful because every valid cutting plane for a given
polyhedron is a $k$-disjunctive cut for some $k$.
\begin{lem}\label{lem:every}
  Let $P \subset \R^{p+q }$ be a polyhedron and $\alpha x + \beta y \leq \gamma$
  be a valid cutting plane. Then $\alpha x + \beta y \leq \gamma$ is a
  $k$-disjunctive cut for some $k \in \N$.
\end{lem}
\begin{proof}
  Let $P$ be a polyhedron and $\alpha x + \beta y \leq \gamma$ be a valid cutting
  plane. The set $M$ that is cut off by the above cutting plane is
  given by
\begin{equation*}
    M:=\{(x,y) \in  P: \alpha x + \beta y > \gamma \}.
  \end{equation*}
  $M$ contains no feasible points of $P_I$. So we have $x \not\in
  \Z^p$ for $(x,y) \in M$. By \autoref{lem:projection} the projection
  of $M$ can be expressed as
  \begin{equation*}
     \prox (M) = \{ x \in \R^p: A^e x \leq b^e, \,  A^l x < b^l\}
  \end{equation*}
  with wlog integral matrices $ A^e, A^l$ with rows $a_i$ satisfying
  $\mathrm{gcd}(a_i) = 1$ and integer vectors $b^e,b^l$. We modify the
  coefficients of the vectors $ b^e, b^l$ by
  \begin{eqnarray*}
    \widetilde b_i^e &=& \lfloor  b_i^e \rfloor + 1, \\
    \widetilde b_i^l &=& \lceil b_i^l \rceil. 
  \end{eqnarray*}
  Altogether we get that $\alpha x + \beta y \leq \gamma$ is a $k$-disjunctive
  cut according to \\
  $ D(k, -(A^e,A^l), -(\widetilde b^e, \widetilde b^l))$.
\end{proof}
It is our goal to compute the mixed integer hull of a given polyhedron
using $k$-disjunctive cuts. Of course this should be done 'as simple
as possible', what means that both the number $k$ of hyperplanes
needed for the disjunctions and the number of iterations in a cutting
plane procedure should be small. At least the latter property can be
easily realized as the next theorem shows.
\begin{thm}\label{thm:2pclos}
  Let $P \subset \R^{p+q}$ be a polyhedron. Then $ P_I =
  P^{(1)}_{2^p}$.
\end{thm}
\begin{proof}
  We will show that every valid inequality $\alpha x + \beta y \leq \gamma$ for
  $P_I$ is a $2^p$-disjunctive cut for $P$. This is sufficient for the
  theorem. Using \autoref{lem:every} we get that $ \alpha x + \beta y \leq \gamma$ is a
  $k$-disjunctive cut with a related disjunction $D(k,d,\delta)$. So
  the claim is shown for $k \leq 2^p$.  Otherwise the number of
  inequalities of the disjunction can be reduced until the required
  limit of $2^p$. Since $D(k,d,\delta)$ is a $k$-disjunction we have
  \begin{eqnarray*}
    \forall x \in \Z^p \ \exists i \in \{1,\ldots,k\}: \ d^ix \leq \delta^i.
  \end{eqnarray*}
  On the other hand we take the set of all integral vectors with the
  property $ d^ix = \delta^i $ for a given $i \in \{1,\ldots,k\}$.
  Either it exists now a vector $\bar{x} \in \Z^p$ with $d^i \bar{x} =
  \delta^i $ and $ d^j \bar{x} > \delta^j, \ \forall j \in
  \{1,\ldots,k\} \setminus \{i\}$, or we can expand the disjunction by
  setting the right hand side of the inequality to $ \delta^i - 1$ and
  repeat this consideration. This may also lead to the case that the
  inequality can be dropped. Therewith we can restrict ourselves to
  disjunctions with the additional condition:
\begin{equation*}
  \forall i \in \{1,\ldots,k\} \ \exists x^i \in \Z^p :
  \ d^ix^i = \delta^i \wedge  d^jx^i > \delta^j, \ 
  \forall  j \in \{1,\ldots,k\} \setminus \{i \} 
\end{equation*}
The set $\conv(\{x^1,\ldots,x^k\})$ contains except for its vertices
$\{x^1,\ldots,x^k\}$ no more integral vectors: Assumed that there was
another integral vector $z \in \conv(\{x^1,\ldots,x^k\})$ we had $ d^i
z \leq \delta^i$ for an $i$. This contradicts the definition of the
vectors $x^i$. So we have constructed a set that contains exactly $k$
integer points as its vertices. This will lead to a contradiction for
$k > 2^p$. Then we had at least two vertices $v,w$ with the additional
property that each component $v_i,w_i, \ i \in \{1,\ldots,n\}$ of both
vectors is either even or odd. So $ \frac{1}{2}(v+w) $ is an integral
vector which is contained in $\conv(\{x^1,\ldots,x^k\})$. This is a
contradiction to the properties of the set.
\end{proof}
We look at an easy example to see that $2^p$-disjunctive cuts are
needed in general to compute the mixed integer hull of a polyhedron in
one step.
\begin{exa}\label{exa:2pclos1}
  We take the $p$-dimensional unit cube $C = [0;1]^p$ and define the
  polyhedron $Q$ by
  \begin{equation*}
    Q = \{x: ax \leq \max_{x \in C} ax, \, a \in \{-1,1\}^p\}.     
  \end{equation*}
  Next we embed $Q$ in the $\R^{p+1}$ and define the polyhedron
  \begin{equation*}
    P = \conv \left\{ \left( \begin{array}{c} x \\ 0 
      \end{array} \right), \frac{1}{2} \mathbf{1} \right\},
  \ x \in Q.
  \end{equation*}
  Of course it is $P_I = C$ and the only valid $k$-disjunction for the
  cutting plane $ x_{p+1} \leq 0$ is defined by the facets of $Q$
  itself.
\end{exa}
As \autoref{thm:2pclos} shows, the mixed integer hull of a general
polyhedron can be 'easily' generated with $2^p$-disjunctive cuts in
theory. Of course for practical issues the use of disjunctions with an
exponential number of defining hyperplane is very expensive. So we
will deal in the following with the second question we mentioned
above, i.e. what kind of $k$-disjunctive cuts we need at least in
computing the mixed integer hull using a repeated application of
k-disjunctive cuts.

\subsection{Approximation property of split cuts}
\label{ssec:convkclos}
Before we further analyze which cuts we need to solve a MILP exactly,
we will deal with the approximation properties of $k$-disjunctive
cuts. We repeat that already using split cuts is sufficient to
approximate the optimal objective function value of any MILP
arbitrarily exact.  Therefor look at the series $ (\gamma^{(i)})_{i
  \in \N}$ of objective function values that is given by
\begin{equation}\label{eq:series}
   \gamma^{(i)} := \max \{cx + hy| \ (x,y) \in P^{(i)}\} 
\end{equation}
for an arbitrary objective function $cx + hy $ that is bounded over
the polyhedron $P$. In detail we get the following 

\begin{thm}\label{thm:finconv}
  Let $P \subset \R^{p+q}$ be a polytope, $cx + hy$ an objective
  function, $ \gamma^* = \max \{cx +hy| \ (x,y) \in P_I \}$ and $
  \gamma^{(i)} $ as defined in \eqref{eq:series}. Then for all
  $\epsilon > 0$ there is an $i_0 \in \N$ with $ | \gamma^{(i_0)} -
  \gamma^*| < \epsilon $.
\end{thm}
\begin{proof}
  A proof of this statement in a slightly different form using a
  repeated variable disjunction can be found in the paper \cite{SM01}
  of Owen and Mehrotra. Moreover the algorithm in this paper also
  gives a constructive proof.
\end{proof}
As we can approximate any optimal objective function value arbitrarily
exact using split cuts, the use of general $k$-disjunctive cuts
becomes necessary for determining exact solutions, only. Moreover we
want to remark that in practical applications already optimizing over
the first split closure often gives a good approximation of the
optimal objective function value. This was in detail investigated by
Balas and Saxena \cite{BS06} for instances from the MIPLIB 3.0 and
several other classes of structured MILP.

\subsection{Solving  MILP exactly}\label{ssec:finite}
We now get back to the question what kind of cuts is needed to solve a
general MILP exactly. As we will see, this depends on the structure of
the projection of the solution space on the $x$-space of integral
variables. For example the important special case of the solution
space being a vertex can be solved just using split cuts. However, we
will see that in general the required number of disjunctive
hyperplanes is exponential in the dimension of the integer space. We
start with the case that the solution set contains relative interior
integer points.
%
\begin{thm}\label{thm:verfin}
  Let $P \subset \R^{p+q}$ be a polyhedron, $cx + hy$ an over $P$
  bounded objective function and $ \gamma^* = \max \{cx +hy| \ (x,y)
  \in P_I \}$. If 
  \begin{equation*}
    \relint( \prox( \{ (x,y) \in P_I: cx + hy = \gamma^*\} )) 
    \cap \Z^p \not= \emptyset
  \end{equation*}
  then there is a $k \in \N$ with $\max \{cx +hy| \ (x,y) \in P^{(k)}
  \} = \gamma^*$.
\end{thm}
\begin{proof}
  If $ \max \{cx +hy| \ (x,y) \in P \} = \gamma^*$ there is nothing to
  show, so let $ \max \{cx +hy| \ (x,y) \in P \} > \gamma^*$. That
  means especially that $ \inte(\prox(M)) \cap \Z^p = \emptyset$ where
  $M := P_I \cap \{ (x,y): cx + hy = \gamma^*\} $ denotes the solution
  set.  Moreover let $x^* \in \relint( \prox(M)) \cap \Z^p$. To proof
  the claim we have to show that $ cx + hy \leq \gamma^*$ is a split
  cut for one of the polyhedra $P^{(k)}, \, k \in \N$.

  Let $A_Ix + G_Iy \leq b_I$ denote these inequalities in the
  representation of the mixed integer hull that constrain the set $M$.
  With \autoref{thm:finconv} we get
  \begin{equation}\label{eq:finconv1}
    \lim_{k \rightarrow \infty} 
    \max \{a_{I,i}x + g_{I,i}y| \ (x,y) \in P^{(k)}\}  = b_{I,i}.
  \end{equation}
  Moreover $x^*$ lies in the boundary of the projection
  $\prox(M^{(k)})$ of each set $ M^{(k)}:= P^{(k)} \cap \{ (x,y): cx +
  hy \geq \gamma^*\}$. As $ M^{(k+1)} \subseteq M^{(k)} $ there exists
  an inequality $px \leq \pi$ that is valid for all of the sets
  $\prox(M^{(k)})$ with the additional property
  \begin{equation}\label{eq:finconv2}
    px = \pi,  \forall x \in \prox(M),
  \end{equation}
  as $x^* \in \relint(\prox(M))$.  If we combine \eqref{eq:finconv1}
  and \eqref{eq:finconv2} we get as direct consequence that $cx + hy
  \leq \gamma^*$ is a split cut to the disjunction $D(p,\pi)$ for some
  $P^{(n)}, n \in \N$.
\end{proof}
After we have seen that split cuts are even sufficient for solving an
important class of MILP exactly, we turn to the general situation. The
idea for finite convergence using $k$-disjunctive cuts in general
consists of the basic principle that there has to exist a
$k$-disjunction $D$ so that the interior of the projection $\prox(M)$
of the solution set is not contained in $D$. If no appropriate
$k$-disjunction exists for all closures $P_k^{(i)}$ then we cannot
achieve a finite algorithm using $k$-disjunctive cuts. On the other
hand, if this condition is satisfied for each face of the solution
set, finite convergence can be shown in the general case.
\begin{thm}\label{thm:mainfin}
  Let $P \subset \R^{p+q}$ be a polyhedron, $cx + hy$ an over $P$
  bounded objective function, $ \gamma^* = \max \{cx +hy| \ (x,y) \in
  P_I \}$ and $M := P_I \cap \{ (x,y): cx + hy = \gamma^*\}$.  If
  there exists for both $M$ and all its faces $f \in F$ with
  $\relint(\prox(f)) \cap \Z^p = \emptyset$ a $k$-disjunction
  $D_f(k,d,\delta)$ with the property
    \begin{equation*}
      x \in \relint(\prox(f)) \ \Longrightarrow x \not\in D_f,
    \end{equation*}
    then there exists a $n \in \N$ with $ \max \{cx + hy| \ (x,y) \in
    P_k^{(n)} \} = \gamma^*$.
\end{thm}
\begin{proof}
  We prove the claim by induction over the dimension $l$ of the
  solution set $M$. We start with $l = 0$. In this case $k = 2 $ can
  always be chosen and the result is a special case of
  \autoref{thm:verfin}. We assume now that the claim is true for
  $l-1,  l \in \N$. \\
  So let $M$ be the solution set of $\max \{cx + hy |\ (x,y) \in P_I
  \}$ and $\dim(M) = l$.  Moreover let for $k \in \N$ exist a
  $k$-disjunction $D(k,d,\delta) $ for $M$ according to the
  assumption. We proof that $cx + hy \leq \gamma^*$ is a k-disjunctive
  cut to the disjunction $D$ for one of the sets $P_k^{(n)}$. Therefor
  we show that it exists a $n \in \N$ such that $ (P_k^{(n)} \cap \{
  (x,y): cx + hy > \gamma^*\}) \cap D = \emptyset$. With the
  disjunction $D$ we define the polyhedron $Q:= \{(x,y): dx \geq
  \delta \wedge cx + hy = \gamma^*\}$. As $cx + hy = \gamma^*$ is a
  supporting hyperplane of $P_I$, there exists $(\widehat x, \widehat
  y)$ with $ c \widehat x + h \widehat y > \gamma^*$ and $ \widehat x
  \not\in D$ such that each of the inequalities $c_f x + h_f y \leq
  \gamma_f$ defined by $(\widehat x, \widehat y)$ and a facet $f$ of
  $Q$ is valid for $P_I$. All inequalities $c_f x + h_f y \leq
  \gamma_f$ at most support $M$ in an under dimensional face. So it
  follows either by induction hypothesis or by \autoref{thm:finconv}
  that the inequalities $c_f x + h_f y \leq \gamma_f$ are valid for
  some $P_k^{(n)}$. Herewith, the condition $ (P_k^{(n)} \cap \{
  (x,y): cx + hy > \gamma^*\}) \cap D = \emptyset$ is satisfied and
  the theorem is proven.
\end{proof}
We remark that it is necessary to involve all the faces of the
solution set in the last theorem, as the following example shows.
\begin{exa}
We define the polyhedron $P \subset \R^{3 + 1}$ through the vertices
\begin{eqnarray*}
  &&(0,0,0,0), (2,0,0,0), (0,2,0,0) \\
  &&(0,0,1,0), (2,0,1,0), (0,2,1,0) \\
  &&\left(\frac{1}{2},\frac{1}{2},\frac{1}{2},\frac{1}{2} \right)
\end{eqnarray*}
For the objective function vector $ (0,0,0,1)$ we get that $ M $ is
contained in a split disjunction, whereas the cut according to the
face $x_3 = 0$ is a $3$-disjunctive cut, only.
\end{exa}
We now deal with the question what cuts we need to solve a general
MILP. Therefor we use special sets that can arise as solution sets of
MILP to give a lower bound of the required number of disjunction
terms. The idea is based on a generalization of \autoref{exa:cks}.
\begin{lem}\label{lem:expneed}
  Let $P \subset \R^{p+q}$ be a polyhedron, $cx + hy$ an over $P$
  bounded objective function, $ \gamma^* = \max \{cx +hy| \ (x,y) \in
  P_I \} > \max \{cx +hy| \ (x,y) \in P \}$ and $M := P_I \cap \{
  (x,y): cx + hy = \gamma^*\}$. If $\prox(M) \subset \R^p$ has $k$
  facets with each containing a relative interior integer point, than
  $ \max \{cx +hy| \ (x,y) \in P_{k-1}^{(i)} \} > \gamma^*, \, \forall
  i \in \N$.
\end{lem}
\begin{proof}
  Let $\prox(M)$ be given with relative interior points $
  \{x^1,\ldots,x^k\} \in \Z^p$ that are contained in pairwise
  different facets. Then we have $ x^{ij} := \frac{1}{2} (x^i + x^j)
  \in \inte(\prox(M))$ for all $i,j, \, i \not=j$. By presumption
  there exists $(x^{ij},y^{ij}) \in P$ with $cx^{ij} + hy^{ij} >
  \gamma^*$. Moreover at least one of the points $(x^{ij},y^{ij})$ is
  not cut off by an arbitrary $k-1$-disjunctive cut. So the cut is
  valid for the set $ Q^{ij} := \conv((x^{ij},y^{ij}), P_I)$. As each
  cut can be classified by this property we get that $ \bigcap_{i
    \not=j} Q^{ij} \subseteq P_{k-1}^{(1)}$. It is clear that
  $\bigcap_{i \not=j} Q^{ij}$ contains a point $(x,y)$ with $cx + hy >
  \gamma^*$. As $P_{k-1}^{(1)}$ satisfies all presumptions and the
  solution set $M$ does not change, the proof follows by induction.
\end{proof}
Therewith we can show now that we need cutting planes to an in the
dimension $p$ exponential number of disjunctive terms to solve a MILP
in general.
\begin{thm}\label{thm:expon}
  Let $P$ be a polyhedron, $cx + hy $ be an over $P$ bounded
  objective function with $ \max\{cx + hy : (x,y) \in P \} =
  \gamma^*$. Then in general
  \begin{equation*}
    \max\{cx + hy : (x,y) \in P_{2^{p-1} + 1 }^{(n)} \} > \gamma^*, 
    \ \forall n \in \N.
  \end{equation*}
\end{thm}
\begin{proof}
  Using \autoref{lem:expneed} it is sufficient to give an integer
  polytope $Q \subset \R^p$ with $p = n + 1$ and at least $2^{n} + 2 $
  facets that contains no interior integer point but in each facet a
  relative interior integer point. Therefor we define $Q$ as the set
  of all $(x,x_{n+1}) \in \R^{n+1}$ with the property:
  \begin{eqnarray*}
    && a x - \pi(a)x_{n+1} \leq 1, \ a \in \{\pm 1\}^n \\
    && 0 \leq x_{n+1} \leq 2
  \end{eqnarray*}
  with $ \pi(a) := |\{i \in \{1,\ldots,n\}:a_i = 1 \} |- 1$. We show
  that $Q$  has the desired properties. \\
  Its vertices are contained in the hyperplane $x_{n+1} = 0$ or
  $x_{n+1} = 2$. For $x_{n+1} = 0$ the related polytope is the
  $n$-dimensional cross polytope. For $x_{n+1} = 2$ the related
  polytope is generated by the vertices $ \mathbf{1} + (n-1)u_i$. The
  last property follows from the fact that $ax \leq 1 $ is active for
  a vector $ \pm u_i$ if, and only if $ ax \leq 1 + 2 \pi(a)$ is
  active for $ \mathbf{1} \pm (n-1)u_i$ by definition of $\pi(a)$.
  So $Q$ is integer. \\
  Let $(z,1) \in \Z^{n+1}$ be given. We take the side constraint $ ax
  \leq 1 + \pi(a)$ with $a_i = 1 \Longleftrightarrow z_i > 0, \ i \in
  \{1,\ldots,n\}$ and get
  \begin{equation*}
    a z = \sum_{1 \leq i \leq n} |z_i| \geq \sum_{z_i > 0} z_i \geq
    |\{ i \in \{1,\ldots,n\}: a_i = 1\}| = 1 + \pi(a). 
  \end{equation*}
  So $(z,1)$ is no interior point of $Q$. Moreover we can see that
  $(z,1) \in Q $ for $ z \in \{0,1\}^n$ and that $(z,1)$ is a relative
  interior point of the facet $ ax - \pi(a) x_{n+1} \leq 1 $ for $a_i
  = 1 \Longleftrightarrow z_i = 1, \ i \in \{1,\ldots,n\}$. As $
  \mathbf{0}$ and $\mathbf{1}$ are relative interior point of $
  x_{n+1} = 0$ and $ x_{n+1} = 2$, $Q$ has all properties.
\end{proof}
\begin{figure}
  \centering 
  \subfigure[$x_3 = 0$ \label{sfig:exa_2p_y0}]
  {\includegraphics[width=37mm]{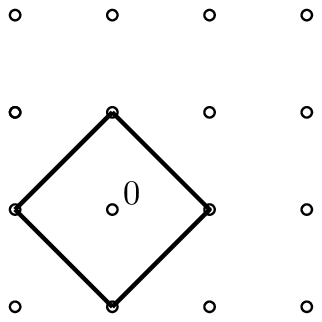}} \hspace*{1.0cm}
  \subfigure[$x_3 = 1$ \label{sfig:exa_2p_y1}]
  {\includegraphics[width=37mm]{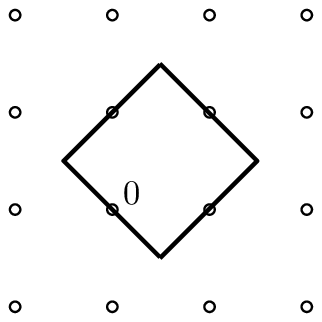}} \hspace*{1.0cm}
  \subfigure[$x_3 = 2$ \label{sfig:exa_2p_y2}]
  {\includegraphics[width=37mm]{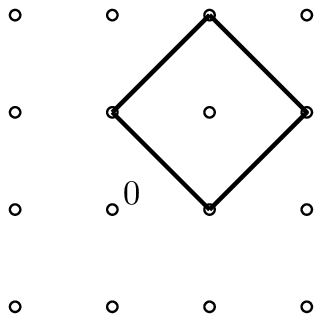}}
  \caption{The set $Q$ constructed in \autoref{thm:expon} for p = 3,
    projected on the $x_3$-space}
  \label{fig:kdis}
\end{figure}
So we have proven that in general at least $2^{p-1} +2$-disjunctive
cuts are required to solve a MILP exactly in a finite number of steps.
We remark that we have an upper bound of $ 2^p$ as shown in
\autoref{thm:2pclos}. With this result we can see that an exact
cutting plane algorithm gets in general very expensive as a large
number of disjunctive hyperplanes has to be computed. Moreover we have
not yet discussed how to determine a cut to the related
$k$-disjunction. As intersection cuts according to basis relaxations
do not generate strong cuts in general, this is an important issue for
practical applications. On the other hand we have seen that a wide
class of problems can even be solved using split cuts. Moreover the
convergence properties of $k$-disjunctive cuts depend on the structure
of the projection of the polyhedron and the objective function. This
fact suggests to use information of the projection in cutting plane
algorithms.

\subsection{Computing $k$-disjunctive cuts}\label{ssec:compkcuts}
At the end of this section we want to give an alternative to compute
strong valid $k$-disjunctive cuts. Unlike the usual generating of
valid cuts for MILP, we need as additional input the vector $(c \ h)$ to
which we want to cut. Moreover we restrict ourselves on shifted
polyhedral cones $ P= \{(x,y): Ax + Gy \leq b\} $ with apex
$(x^*,y^*), x^* \not\in \Z^p$ and assume that the function $ cx
+ hy$ attains its unique maximum at $(x^*,y^*)$ with value $\gamma^*$.
In this situation we can describe $(x^*,y^*)$ as polyhedron given by
the (over-determined) system
\begin{equation}\label{eq:compkcut1}
P_{\gamma^*} := \left\{  
\left( \begin{array}{c}   A \\ -c \end{array}\right)x +
  \left( \begin{array}{c} G \\ -h \end{array}\right)y \leq 
   \left( \begin{array}{c} b \\ -\gamma^* \end{array}\right)
\right\}.
\end{equation}
Using \autoref{lem:projection}, the projection of the above system on the
$x$-space - that is equal to $x^*$ -  is given by
  \begin{equation}\label{eq:compkcut2}
    \prox(P_{\gamma^*}) = \left\{ x \in \R^p: v^r 
    \left( \begin{array}{c} A \\  -c \end{array} \right)
    x \leq v^r 
    \left( \begin{array}{c} b  \\ -\gamma^* \end{array} \right)
    , \ \forall r \in R \right\} 
  \end{equation}
with $R$ being the set of extreme rays of the cone 
\begin{equation}\label{eq:compkcut3}
  Q = \left\{ v \in \R^{m+1}: 
     \left( \begin{array}{c} G  \\ -h \end{array} \right)^T v = 0, \ 
     v \geq 0 \right\}.
\end{equation}
As the cone $Q$ is rational, we can assume that the extremal rays
$v^r$ are elements of the additive group
\begin{equation}\label{eq:compkcut4}
  \mathcal{G}_{P_{\gamma^*}} := \left\{w \in \Q^{m+1}: 
    w \left( \begin{array}{c} A \\ -c \end{array} \right)  \in \Z^{m+1} 
  \right\}.
\end{equation}
Therewith we can use the above polyhedral description of
$\prox(P_{\gamma^*})$ to define a valid $k$-disjunction for $P$, that
does not contain the apex $(x^*,y^*)$. We do this by rounding up the
right hand sides of the defining constraints of $\prox(P_{\gamma^*})$
in \eqref{eq:compkcut2}.
\begin{lem}\label{lem:Rdis}
  Let $P,(x^*,y^*),(c \ h),\gamma^*$ and $ \prox(P_{\gamma^*})$ as
  defined above. Moreover define for $r \in R$
  \begin{eqnarray*}
    d^r &:=&  v_{1,\ldots,m}^r A   - v_{m+1}^r c \\
    \delta^r &:=& \lfloor v_{1,\ldots,m}^r b  - v_{m+1}^r \gamma^* \rfloor + 1
  \end{eqnarray*}
  with $ v^r = (v_{1,\ldots,m}^r, v_{m+1}^r)$.  Then $
  D(|R|,-d,-\delta)$ is a valid $|R|$-disjunction for $P$ that does
  not contain $(x^*,y^*)$.
\end{lem}
\begin{proof}
  By definition it is $\max\{ cx +hy: (x,y) \in P_I\} < \gamma^*$.  So
  there is an $\epsilon > 0$ such that $cx + hy \leq \gamma^* -
  \epsilon $ is valid but not optimal for $P_I$. Moreover the
  inequalities defining $ \prox(P_{\gamma^* - \epsilon})$ are given by
  $ v_{1,\ldots,m}^r A - v_{m+1}^r c $ with right hand sides
  $v_{1,\ldots,m}^r b - v_{m+1}^r (\gamma^* - \epsilon)$. The set $
  \prox(P_{\gamma^* - \epsilon})$ contains no integer points, so for
  all $x \in \Z^p$ there is a $r \in R$ with 
  \begin{equation*}
    d^r x >  v_{1,\ldots,m}^r b  - v_{m+1}^r \gamma^*.
  \end{equation*}
  Therewith is follows that the polyhedron $ \{x \in \R^p: dx \leq
  \delta\} $ contains no integer point in its interior. This is
  equivalent to the set $ D(|R|,-d,-\delta)$ being a valid
  $|R|$-disjunction. Moreover, as the right hand side of each defining
  hyperplane of the projection has been enlarged by the definition of
  $D$ it is obvious that $(x^*,y^*$) is not contained in the
  disjunction. This proofs the lemma. We remark that the above
  definitions and the proof is similar to \autoref{lem:every}.
\end{proof}
So we have found a $k$-disjunction that can be used to cut off the
current LP solution $(x^*,y^*)$. As we have mentioned at the beginning
of this section we want to cut to the vector $(c \ h)$. We can do this
now using the right hand side $\delta$ of the disjunction $D$. As for
$ v_{m+1}^r > 0$ the value of $\delta^r$ depends on the objective
function value, we can compute the objective function value that
corresponds to the value of $\delta^r$ that we have got by the
rounding operation. The inequalities of $D$ whose right hand sides
$\delta^r$ are independent of the value of $\gamma$ can be omitted in
this considerations. Taking the maximum of the related objective
function values for all constraints gives us a valid cut to the vector
$(c \ h)$. As the disjunction does not contain $x^*$, we can ensure
that the current solution is cut off.
\begin{thm}\label{thm:Rdiscut}
  Let $P,(x^*,y^*),(c \ h),\gamma^*, \prox(P_{\gamma^*})$ and $
  D(|R|,-d,-\delta)$ as defined above. Let for $r \in R$ with
  $v_{m+1}^r > 0$
  \begin{equation*}
    \gamma^r := \frac{\delta^r - v_{1,\ldots,m}^r b }{ - v_{m+1}^r}
   \end{equation*}
   and $\widehat \gamma = \max\{\gamma^r: r \in R, \, v_{m+1}^r >
   0\}$. Then $cx + hy \leq \widehat \gamma $ is valid for $P_I$ and
   $cx^* + hy^* > \widehat \gamma$.
\end{thm}
\begin{proof}
  The validity of the inequality $cx + hy \leq \widehat \gamma $
  follows directly using \autoref{lem:Rdis}, as it is a disjunctive
  cut according to $ D(|R|,-d,-\delta)$ by definition. Equally it
  follows that $cx^* + hy^* > \widehat \gamma$.
\end{proof} 
As we have finished the derivation of the $k$-disjunctive cut, we want
to add some remarks. Using the projection as $k$-disjunction, we solve
the problem how to find a suitable $k$-disjunction for cutting in
general. This relates both to the selection of the number $k$ and the
selection of the defining hyperplanes of the disjunction. Moreover we
have seen in the last subsections, that using information of the
projection can be useful. On the other hand, the projection that we
use corresponds to the predisposed cutting vector. So the selection of
a suitable $k$-disjunction is partially shifted to the selection of
the cutting vector. Here it is i.e. open how to choose cutting vectors
to get deep cuts in general.  However, for solving a given MILP we
will see in the next section that this approach leads to a finite
algorithm if we use the objective function vector.

\section{Algorithm}\label{sec:algo}
We now turn to an algorithmic application of the previous results and
want to present an exact algorithm that solves a bounded MILP in
finite time. It is based on a series of mixed integer Gomory cuts that
is mixed with certain $k$-disjunctive cuts which are required as
discussed in \autoref{ssec:finite}. The $k$-disjunctive cuts we use
here are similar to the ones we introduced in
\autoref{ssec:compkcuts}, using the objective function as the vector
to which we cut. As the assumptions that we have made there for the
$k$-disjunctive cuts are in general not satisfied, we have to do some
modifications. So we will define $k$-disjunctive cuts over general
polyhedra $P$ for an arbitrary cut vector $(c \ h)$.
We discuss the details of the generalization. It is clear that the
equalities \eqref{eq:compkcut2}, \eqref{eq:compkcut3},
\eqref{eq:compkcut4} also describe the projection for $P$ being a
general polyhedron and $\gamma^*$ being an arbitrary value of the
objective function. Even the derivation of a valid $k$-disjunction and
a valid $k$-disjunctive cut, respectively, is true, if the value
$\gamma^*$ of the objective function $cx +hy$ is not optimal for
$P_I$. However, for the application in the algorithm we will define a
slightly weaker version of the $|R|$-disjunctive cut that does not
always cut off the current LP solution, but can be used more general.
We do this in the next
\begin{thm}\label{thm:genRdiscut}
  Let $P$ be a polyhedron and $\gamma^*$ such that $c x + h y \leq
  \gamma^*$ is valid for $P_I$ but not for $P$. Define the
  $|R|$-disjunction $ D(|R|,-d,-\delta)$ using the equalities
  \eqref{eq:compkcut1}, \eqref{eq:compkcut2}, \eqref{eq:compkcut3},
  \eqref{eq:compkcut4} with
 \begin{eqnarray*}
    d^r &:=&  v_{1,\ldots,m}^r A   - v_{m+1}^r c \\
    \delta^r &:=& \lceil v_{1,\ldots,m}^r b  - v_{m+1}^r \gamma^* \rceil
  \end{eqnarray*}
  and let $\widehat \gamma = \max\{\gamma^r: r \in R, \, v_{m+1}^r >
  0\}$ analog to \autoref{thm:Rdiscut}. Then $c x + h y \leq \widehat
  \gamma $ is a valid cutting plane for $P_I$. 
\end{thm}
\begin{proof}
  By assumption $\prox(P_{\gamma^*})$ cannot contain an integral point
  in its interior. So rounding up the right hand sides gives a valid
  $|R|$-disjunction. Therefor $c x + h y \leq \widehat \gamma $ is a
  valid cutting plane analog to the proof of \autoref{thm:Rdiscut}.
\end{proof}
We go on with the single steps of the algorithm. We start with the
usual mixed integer Gomory algorithm as long as we get either a
feasible solution of the MILP or a solution of the LP relaxation that
has a lower objective function value. This happens in finite time by
\autoref{cor:gomory} if we restrict ourselves to polytopes. If the
objective function value has decreased we can apply
\autoref{thm:genRdiscut} and compute a valid $|R|$-disjunctive cut
using the objective function as vector to which we cut. Now we can
apply the Gomory algorithm to the modified program again, until either
a feasible solution is found or the objective function value
decreases, and use \autoref{thm:genRdiscut} again. In this way we get
an algorithm that finitely terminates with an optimal solution to the
given MILP or detects infeasibility. The formal algorithm is stated in
\autoref{alg:exact}.
 \begin{algorithm}\caption{Exact cutting plane algorithm}  
 \label{alg:exact}
   \begin{algorithmic}[1]
     \State \textbf{Input:} bounded MILP \eqref{eq:MILP} 
     \State \textbf{Output:} "optimal solution $(x^*,y^*)$" or 
     "problem infeasible" if no solution exists; 
     \State 
     \State $(x^*, y^*) := \argmax\{cx + hy: (x,y) \in P\}$; 
     \State $ \gamma^* := \max\{cx + hy: (x,y) \in P\}$;
     \State
     \If{$ P = \emptyset$} 
     \State "problem infeasible"; \textbf{break}
     \EndIf
     \If{$ x^* \in \Z^p$} \State "optimal solution 
     $ (x^*,  y^* )$"; \textbf{break}
     \EndIf
     \State
     \While{$x^* \not\in \Z^p$}
     \State $ \gamma := \gamma^*$; 
     \While{$\gamma^* = \gamma$} 
     \State Compute Gomory cut $ \alpha^1 x + \alpha^2 y \leq \beta$ 
     to $P, (x^*, y^*)$ by least index rule;
     \State $ P := P \cap \{(x,y): \alpha^1 x + \alpha^2 y \leq \beta\}$; 
     \State $( x^*,  y^*) := \argmax\{cx + hy: (x,y) \in P\}$; 
     \State $\gamma^* := \max\{cx + hy: (x,y) \in P\}$;
      \If{$ P  = \emptyset$} 
     \State "problem infeasible"; \textbf{break}
     \EndIf
     \If{$ x^* \in \Z^p$} \State "optimal solution 
     $( x^*, y^*)$"; \textbf{break}
     \EndIf
     \EndWhile
     \State 
     \State Compute $\widehat \gamma = \max\{\gamma^r : r \in
     R, v_{m+1}^r > 0 \}$ according to \autoref{thm:genRdiscut}
     for $P, (c \ h ), \gamma^*$;
     \State $\gamma^* = \widetilde \gamma$;
     \EndWhile
     \State
 \end{algorithmic}
 \end{algorithm}
\begin{thm}\label{thm:finalg}
  Let a bounded MILP \eqref{eq:MILP} be given. Then
  \autoref{alg:exact} either finds an optimal solution or detects
  infeasibility in a finite number of steps.
\end{thm}
\begin{proof}
  The proof follows immediately with the following two facts: Every
  while loop (16) to (27) has only finite many iterations by
  \autoref{cor:gomory} as $P$ is bounded by presumption. Similarly the
  outer while loop (14) to (31) has only finite many iterations as the
  possible number of different values $\widehat \gamma$ is finite. 
\end{proof}
Before we further discuss the algorithm we will give two examples. We
start with repeating \autoref{exa:cks}:
\begin{exa}
  Let again the MILP
  \begin{eqnarray*}
    && \max y \\
    && -x_1 + y \leq 0 \\
    && -x_2 + y \leq 0 \\
    &&  x_1 + x_2 + y \leq 2 \\
    && x_1, x_2 \in \Z
  \end{eqnarray*}
  with the optimal solution $(\frac{2}{3},\frac{2}{3},\frac{2}{3})$ of
  the LP relaxation be given. The mixed integer Gomory cuts according
  to $x_1$ and $x_2$ are given by $ -x_1 + 2 y \leq 0 $ and $ -x_2 + 2
  y \leq 0 $ with the new LP solution $ (\frac{4}{5},\frac{4}{5},
  \frac{2}{5})$. As the value of the objective function has decreased,
  we compute $\widetilde \gamma$ as in \autoref{thm:genRdiscut}. The
  extremal rays of the cone $\{(\begin{array}{cccc} 1 & 1 & 1 & -1
    \end{array})y = 0, \  y \geq 0\}$
  are the three vectors
  \begin{equation*}
    (1,0,0,1),(0,1,0,1),(0,0,1,1),
  \end{equation*}
  so the projection $\prox (P_\gamma)$ is given by
  \begin{eqnarray*}
    -x_1 &\leq& 0 - \gamma\\
    -x_2  &\leq& 0 - \gamma \\
    x_1 + x_2  &\leq& 2 - \gamma 
 \end{eqnarray*}
 Inserting the current value $ \gamma = \frac{2}{5}$ of the objective
   function and rounding gives \\ $\max_{r =1,2,3}  \tilde \gamma^r =
   \max\{0,0,0\} = 0$. After applying the related cut $ y \leq 0$ we
   get as new LP solution the feasible point $(2,0,0)$ and the
   algorithm stops with an optimal solution. 
\end{exa}
Second we show how the algorithm works for the example of Owen and
Mehrotra \cite{SM01}. For this ILP the usual mixed integer Gomory
algorithm does not converge to the optimum.
\begin{exa}
  Let the ILP
 \begin{eqnarray*}
     \max&&x_1 + x_2 \\
    && 8x_1 + 12x_2 \leq 27 \\
    && 8x_1 + 3 x_2 \leq 18 \\
    && x_1, x_2 \geq 0 \\
    && x_1, x_2 \in \Z
  \end{eqnarray*}
  with the initial LP solution $(\frac{15}{8}, 1) $ be given. After
  applying the first possible cut to $x_1$ the value of the objective
  function decreases and we can go to the second step of the
  algorithm. As we have an ILP it is $ \prox(P_{\gamma}) = P_{\gamma}
  $ with $ P_{\gamma} $ given by
  \begin{eqnarray*}
     - x_1 - x_2 &\leq& -\gamma \\
     8x_1 + 12x_2&\leq& 27 \\
     8x_1 + 3 x_2 &\leq& 18 \\
     x_1, x_2 &\geq& 0 \\
  \end{eqnarray*}
  By rounding we get finally the valid cut $ x_1 + x_2 \leq 2 $ that
  relates to the optimal objective function value.

  The result of this example is typical for applying the algorithm on
  ILP. In this case we have to presume that all input data is integral
  and the $k$-disjunctive cut to the objective function reduces to the
  Chv\'atal Gomory cut of the objective function vector.
\end{exa}
Concluding we want to discuss the algorithm. We have seen in the last
example that for an ILP the $k$-disjunctive cut reduces to an integer
Gomory cut to the objective function. So the whole algorithm can be
seen as a variant of the pure integer Gomory algorithm in this case.
The crucial fact for finite convergence of the integer algorithm is
the possibility to add cuts both to the objective function and to each
variable if they are not integral. Using $k$-disjunctive cuts to the
objective function we have now the possibility to add cuts to the
objective function in the case of MILP as well. Therewith we obtain a
convergent algorithm in analogy to the integer case.

Of course the complex part of the algorithm consists in computing the
$k$-disjunctive cut as the number $|R|$ of extreme rays $v^r$ of the
cone $Q$ grows exponentially. So an efficient algorithm for computing
the extreme rays of the related cone is required. Moreover we have to
ensure that the computed rays satisfy the integrality constraints,
i.e. are contained in the group $ \mathcal G_{P_{\gamma^*}}$. Therefor
we can presume in practical applications the coefficients of the
matrix $A$ and the vector $c$ to be integer. Then the integrality
constraints are satisfied, if all of the extreme rays $v^r$ are
integer. However, we will not further deal with this issue here, but
refer to the papers of Henk and Weismantel \cite{HW96} and of Hemmecke
\cite{Hem06} and the references therein. They state several algorithms
for this and the similar problem of computing Hilbert bases of
polyhedral cones.

At last we want to give a further interpretation of the algorithm.
Therefor we assume that the feasible domain $P$ is full dimensional
and bounded.  One can see that in this case we can always choose an
optimal solution of the MILP such that $q$ defining inequalities of
$P$ are active.
So the solution is contained in a $ (p+q) - q = p $-dimensional face
of $P$. Therefor we can solve the MILP by solving each of the related
$p$-dimensional subproblems and taking the best solution. Moreover the
set of feasible solutions in each $p$-dimensional face is discrete in
general, so solving a MILP for a $p$-dimensional face can be
interpreted as solving an ILP, as we could apply a suitable affine
transformation. This means that solving a MILP can be seen as parallel
solving of several ILP. Especially every valid cutting plane for $P_I$
is even valid for each of the discrete subproblems. Therefor we need
information of the related discrete subproblems if we want to generate
strong valid cuts.  As the number of $p$-dimensional faces of $P$
grows exponentially, this interpretation also gives another reasoning
that we need $k$-disjunctive cuts with an exponential number of
defining disjunctive hyperplanes to solve general MILP. Within the
algorithm we can find the $p$-dimensional subproblems in the facets of
the polyhedron $\prox(P_{\gamma^*})$, where the value of the right
hand side $\delta^i$ can be related to the current objective function
value of the subproblem.

\bibliographystyle{plain} 
\bibliography{Literatur}
\end{document}